\documentclass{amsart}
\usepackage{amsthm,amsmath,amssymb}
\newtheorem{Def}{Definition}[section]
\newtheorem{Lem}{Lemma}[section]
\newtheorem{Thm}{Theorem}[section]

\newtheorem{Ex}{Example}[section]

\PassOptionsToPackage{linktocpage}{hyperref}
\usepackage{hyperref}
\numberwithin{equation}{section}
\begin{document}

\title{On the Trajectories of the $3x+1$ Problem}

\author{Roy Burson}
\address{Department of Mathematics, California State University Northridge, California, 91330}
\email{roy.burson.618@my.csun.edu.}

\date{May 17, 2020}



\begin{abstract}
This paper studies certain trajectories of the Collatz function. I show that if $n\sim 3n+2$ for each odd number $n$ then every positive integer $n \in \mathbb{N}\setminus 2^{\mathbb{N}}$ has the representation $$n=\left(2^{a_{k+1}}-\sum_{i=0}^{k}{2^{a_i}3^{k-i}}\right)/ 3^{k+1}$$ where $0\le a_0 \le a_1 \le \cdot \cdot \cdot \le a_{k+1}$.  As a consequence, in order to prove Collatz Conjecture it is sufficient to prove $n\sim 3n+2$ for each $n\in \mathbb{N}\setminus 2^{\mathbb{N}} $. This is the main result of the paper.
\end{abstract}

\maketitle

\section{Introduction}\label{S: Introduction}
\indent Some problems in mathematics are easy to state but take very complex tools to prove and often it takes new tools to be developed. The \textit{Collatz Conjecture} is either one of this type or else it might be a conspiracy theory to slow down the field of mathematics (a joke that spread across Yale University). The \textit{Collatz Conjecture} is a well known unsolved mathematical problem that concerns the recursive behavior of the function

$$ T(n)=
\begin{cases}
\frac{n}{2} & \textit{if $n\equiv 0(\bmod{2})$}\\
3n+1 & \textit{if $n\equiv 1(\bmod{2})$}\\
\end{cases}$$

over the set of integers $\mathbb{Z}$. In this paper I focus on the Collatz function over the positive whole numbers $\mathbb{N}$. The problem is most commonly referred to as the \textit{"$3x+1$" problem}. The history and exact origin of the problem is somewhat vague. Some early history on the problem is discussed by Jeffrey C. Lagarias at \url{http://www.cecm.sfu.ca/organics/papers/lagarias/paper/html/node1.html}. Lagarias gives 197 different documentations on the topic in his annotated bibliographies \cite{arxiv/0309224v13} and \cite{arxiv/0608208v6}. Dr. Lothar Collatz is credited for the discovery of the problem during his career as a student. Dr. Lothar Collatz even asserts himself that he was the first to study this problem in his letter \url{http://www.cecm.sfu.ca/organics/papers/lagarias/paper/html/letter.html}.\\ 

\indent The problem states that $\forall n\in \mathbb{N}$ there is a value $k\in \mathbb{N}$ such that $T^k(n) = 1$ were $T^k$ is the $kth$ application of the map $T$. That is $T^k$ is the map $T^k: \mathbb{N}\times\mathbb{N}\rightarrow \mathbb{N}$ defined by the rule $$(n,k)\mapsto(\underbrace{T\circ T\circ T \cdots \circ T}_{k-times}) \circ (n)$$ 

\begin{Ex}\label{Example 1}
For the value $n=11$ we can view the iterations with arrows (to indicate direction) as followed: $$11\rightarrow 34 \rightarrow 17 \rightarrow 52 \rightarrow 26 \rightarrow 13 \rightarrow 40 \rightarrow 20 \rightarrow 10 \rightarrow 5 \rightarrow 16 \rightarrow \cdot \cdot \cdot \rightarrow 1$$ We can also visualize this backwards by reversing the operations of the map and traces our steps in the reverse direction. Doing so for this example we have the following: \begin{align*}
1 & \rightarrow \cdot \cdot \cdot \rightarrow \frac{2^4-1}{3}\rightarrow  \frac{2^5-2}{3}\rightarrow \frac{2^6-2^2}{3}\rightarrow \frac{2^7-2^3}{3} \rightarrow \frac{2^7-2^3-3}{3^2}\\
& \rightarrow \frac{2^8-2^4-2\cdot3}{3^2}\rightarrow \frac{2^9-2^5-2^2\cdot3}{3^2}\rightarrow \frac{2^9-2^5-2^2\cdot3-3^2}{3^3}\\
& \rightarrow \frac{2^10-2^6-2^3\cdot3-2\cdot 3^2}{3^3} \rightarrow \frac{2^{10}-2^6-2^3\cdot3-2\cdot 3^2-3^3}{3^4}=11
\end{align*} 

\end{Ex} 
It has been verified that all natural numbers $n< 87 \cdot 2^{60}$ iterate to the value $1$ under the Collatz function. A neat discussion about the empirical results and the record-holders are discussed by Tom\'as Oliveria e Silvia at his home page \url{http://sweet.ua.pt/tos/3x+1.html}, and by Eric Roosendaal \url{http://www.ericr.nl/wondrous/index.html}. A complete list of the record holders can be found here \url{http://www.ericr.nl/wondrous/pathrecs.html} which was accomplished by the yoyo@home project in 2017. In \cite{arxiv/0309224v13} and \cite{arxiv/0608208v6} Lagarias also mentioned the empirical evidence of the problem given by Oliveir\'a e Silva. An interesting study on properties of divergent trajectories of the Collatz function (if one exist) can be found here \url{http://www.csun.edu/~vcmth02i/Collatz.pdf}. A similar result of this paper is given by Charles C. Cadogan in his works \cite{MC44} and \cite{MC11}, which can also be found in Lagarias's bibliography \cite{arxiv/0608208v6}. Some more interesting discussions and findings that relate to this work can be found in \cite{Qufu12, MS63, AMJC, AMM, AMSJC}.\\

\section{Terminology}\label{S: Terminology}
In order to prove the main result of the paper the following definitions are needed. 
\begin{Def}\label{D: Representation}
let $m\in \mathbb{N}$. Write $m\equiv \mathcal{R}$ if and only if the number $m$ has the representation $$m =\left(2^{a_{k+1}}-\sum_{i=0}^{k}{2^{a_i}3^{k-i}}\right)/ 3^{k+1}$$ where $\left(a_i\right)$ is a monotonically increasing sequence of positive integers for $i\le k+1$. 
\end{Def}

\begin{Def}\label{D: Forward Orbit}
The \textit{Trajectory} or \textit{Forward Orbit} of a positive integer $n$ is the set $$O^+(n)=\{n, T(n), T^2(n), \cdot \cdot \cdot \}$$ were $T:\mathbb{N}\rightarrow \mathbb{N}$ is the Collatz function defined by $n\rightarrow \frac{n}{2}$ if n is even and $n\rightarrow 3n+1$ if n is odd, and $T^k$ is the function $T: \mathbb{N}\times \mathbb{N} \rightarrow \mathbb{N}$ defined by $$(n,k)\mapsto \underbrace{(T\circ T \circ \cdot \cdot \cdot \circ T)}_{k-times}\circ (n)$$ 
\end{Def}

\begin{Def}\label{D:Coalesce}
Given two integers $n_1$ and $n_2$ define the relation $n_1\sim n_2$ if and only if the two trajectories $O^+(n_1)$ and $O^+(n_2)$ coalesce, i.e. $O^+(n_1) \cap O^+(n_2)\neq \emptyset$.
\end{Def}

\section{ Backwards Iterations and Integer Representations} \label{S3: Backwards Iterations }
As shown by \cite{M43,M37, B.ohm, Cran, M148}, and discussed in \cite{Goodwin} we have that 

$$ \forall m \in \mathbb{N}\setminus 2^{\mathbb{N}}, 1\in O^+(m) \Longleftrightarrow m=\left(2^{a_{k+1}}-\sum_{i=0}^{k}{2^{a_i}3^{k-i}}\right)/3^{k+1}$$ 

for some sequence $(a_j)$ where $0\le a_0 \le a_1 \le a_2 \le \cdot \cdot \cdot a_{k+1}$. If $m$ has this representation given above then $m$ cannot be a power a $2$ because any element of the set $2^{\mathbb{N}}=\{2,2^2,2^3,\cdot\cdot \cdot\}$ cannot be represented in this form. This can easily be verified by observing that the numbers that can be written in the form $\left(2^{a_{k+1}}-\sum_{i=0}^{k}{2^{a_i}3^{k-i}}\right)/3^{k+1}$ are exactly those numbers obtained by iterating the functions $n\mapsto 2n$ and $n\mapsto \frac{n-1}{3}$ in succession, beginning at the number $n=1$, were as each function is applied at least once and the function $n\mapsto \frac{n-1}{3}$ is never applied twice in succession. The numbers in the set $2^{\mathbb{N}}$ can only be obtained by iterating the function $n\rightarrow 2n$, beginning at $n=1$, so that the function $n\rightarrow \frac{n-1}{3}$ is never applied. This shows that no integer in the set $2^{\mathbb{N}}$ has the representation $\left(2^{a_{k+1}}-\sum_{i=0}^{k}{2^{a_i}3^{k-i}}\right)/3^{k+1}$. The next lemma shows that if a number $n$ has this representation then the number $2n$ and $\frac{2n-1}{3}$ must also have the same representation.

\begin{Lem}\label{L: Lemma1} (Closure of $\mathcal{R}$)
If $n\equiv \mathcal{R}$ then $2n\equiv\mathcal{R}$ and if  $n\equiv \mathcal{R}$ then $\frac{2n-1}{3}\equiv \mathcal{R}$. 
\begin{proof} (Direct proof)
First suppose $n\equiv \mathcal{R}$. Then there is a sequence $(a_j)$ so that $$n=\left(2^{a_{k+1}}-\sum_{i=0}^{k}{2^{a_i}3^{k-i}}\right)/3^{k+1}$$ where $0\le a_0 \le a_1 \le a_2 \le \cdot \cdot \cdot a_{k+1}$. Write 

$$2n = 2\left(\left(2^{a_{k+1}}-\sum_{i=0}^{k}{2^{a_i}3^{k-i}}\right)/3^{k+1}\right)= \left(2^{a_{k+1}+1}-\sum_{i=0}^{k}{2^{a_i+1}3^{k-i}}\right)/3^{k+1}$$

Define the sequence $(b_j)$ by $b_j=a_{j+1}+1$ for each $0\le j\le k$ so that $(b_j)$ is also positive and increasing. Then we have $$2n = \left(2^{b_k+1}-\sum_{i=0}^{k}{2^{b_j}3^{k-i}}\right)/3^{k+1}\equiv \mathcal{R}$$ Therefore $2n\equiv \mathcal{R}$. Now suppose that $n\equiv\mathcal{R}$ then by the above $2n\equiv \mathcal{R}$. Write $$2n=\left(2^{a_{k+1}}-\sum_{i=0}^{k}{2^{a_i}3^{k-i}}\right)/3^{k+1}$$ were as $0\le a_0\le a_1\le \cdot \cdot \cdot a_{k+1}$. Then $$\frac{2n-1}{3} =\frac{\left(2^{a_{k+1}}-\sum_{i=0}^{k}{2^{a_i}3^{k-i}}\right)/3^{k+1}-1}{3}$$ Define the sequence $(b_j)$ by 
$$
b_j = 
\begin{cases}
0& \text{if $j=0$}\\
a_{j-1} & \text{if $1\le j\le k+1$}
\end{cases}
$$

Then $$\frac{2n-1}{3} = \left(2^{b_{k+1}}-\sum_{i=0}^{k+1}{2^{b_i}3^{k+1-i}}\right)/3^{k+2}$$
and the sequence $(b_j)$ remains monotonically increasing. Therefore, it follows that $\frac{2n-1}{3}\equiv \mathcal{R}$
\end{proof}
\end{Lem}

\begin{Lem} \label{L: Lemma2}
Let $a\in \mathbb{N}$. Then $$3^{\frac{a}{2}+1}+2<2^a+1$$ for all $a\ge 8$.
\begin{proof}(By induction)
The proof follows by induction. Let $$S=\{1,2,3,\cdot \cdot \cdot, 7\}\cup \{a\in \mathbb{N}: 3^{\frac{a}{2}+1}+2\le 2^a+1 \}$$ The value $a=8$ is in $S$ since $$245 = 3^{5}+2 < 2^{8}+1 = 257$$ Now if $a\ge 8 \in S$ then \begin{align*}
3^{\frac{a}{2}+1}+2\le 2^a+1 & \Rightarrow 3^{\frac{1}{2}}(3^{\frac{a}{2}+1}+2)\le 2(2^a+1)\\
& \Rightarrow 3^{\frac{1}{2}}(3^{{\frac{1}{2}}^{(a+2)}}+2)\le 2(2^a+1)\\
& \Rightarrow 3^{{\frac{1}{2}}^{(a+3)}}+23^{\frac{1}{2}}\le 2^{a+1}+2\\
& \Rightarrow 3^{{\frac{1}{2}}^{(a+3)}}+3\le 2^{a+1}+2\\
& \Rightarrow 3^{{\frac{1}{2}}^{(a+3)}}+2\le 2^{a+1}+1\\
& \Rightarrow 3^{\frac{a+1}{2}+1}+2\le 2^{a+1}+1
\end{align*}
Therefore $a+1\in S$ so that $S=\mathbb{N}$.
\end{proof}
\end{Lem}

\begin{Lem} \label{L: Lemma3}
Let $n\in 2\mathbb{N}=\{2n : n \in \mathbb{N}\}$ and write $n$ is its conical form $$n = 2^{\epsilon}\prod_{i=1}^{k}{p_{\beta_i}^{\alpha_i}}$$ Then there exist a value $k\in \mathbb{N}$ so that $$T^k(n+1) = \begin{cases}
3^{\frac{\epsilon}{2}}\left(\prod_{i=1}^{k}{p_{\beta_i}^{\alpha_i}}\right)+1 & \textit{if $\epsilon \equiv 0(\bmod{2}$})\\
3^{\lfloor\frac{\epsilon}{2}\rfloor +1}\left(\prod_{i=1}^{k}{p_{\beta_i}^{\alpha_i}}\right)+2 & \textit{if $\epsilon \equiv 1(\bmod{2})$}\\
\end{cases}$$
\begin{proof}
Let $n\in 2\mathbb{N}$ and write $n = 2^{\epsilon}\prod_{i=1}^{k}{p_{\beta_i}^{\alpha_i}}$ for some $\epsilon \in \mathbb{N}$. First, assume $\epsilon \equiv 0(\bmod{2})$ and take $k= 3\frac{\epsilon}{2}$. Then the claim is that $T^k(n+1)= 3^{\frac{\epsilon}{2}}\left(\prod_{i=1}^{k}{p_{\beta_i}^{\alpha_i}}\right)+1$. Write $n+1 = 2^{\epsilon}\prod_{i=1}^{k}{p_{\beta_i}^{\alpha_i}}+1$. Define the sequence $(a_i)$ by the rule 
$$a_{i+1} = \begin{cases}
a_{i}+1 & \text{if $i\equiv 0(\bmod{2})$}\\
a_{i}+2 & \text{if $i\equiv 1(\bmod{2})$}\\
\end{cases}
$$ 
were as $a_1 = 1$. Now since $n+1$ must be odd by computing $T^k$ consecutively we find
\begin{align*}
n+1 =  2^{\epsilon}\prod_{i=1}^{k}{p_{\beta_i}^{\alpha_i}}+1  \Rightarrow T(n+1) &=  32^{\epsilon}\prod_{i=1}^{k}{p_{\beta_i}^{\alpha_i}}+2^2\\
 \Rightarrow T^3(n+1) &=  3^22^{\epsilon-2}\prod_{i=1}^{k}{p_{\beta_i}^{\alpha_i}}+1\\
 \Rightarrow T^4(n+1) &=  3^32^{\epsilon-2}\prod_{i=1}^{k}{p_{\beta_i}^{\alpha_i}}+2^2\\
\Rightarrow T^{6}(n+1) &=  3^42^{\epsilon-4}\prod_{i=1}^{k}{p_{\beta_i}^{\alpha_i}}+1\\ 
& \indent \vdots \\
\Rightarrow T^{k}(n+1) &=  3^i\prod_{i=1}^{k}{p_{\beta_i}^{\alpha_i}}+1\\ 
\end{align*}
for some $i$ and $k$ because $\epsilon$ is even. Actually we know $i =\frac{\epsilon}{2}$ and $k=a_i$ (were $(a_i)$ was defined above) so $$T^{k}(n+1) = 3^{\frac{\epsilon}{2}}\prod_{i=1}^{k}{p_{\beta_i}^{\alpha_i}}+1$$ as desired. Now if  $\epsilon\equiv 1(\bmod{2})$ then we may write $\epsilon = \epsilon^{\prime} +1$ were $\epsilon \equiv 0(\bmod{2})$. Hence, we see that \begin{align*}
T^k(n+1) &=T^k(2^{\epsilon}\prod_{i=1}^{k}{p_{\beta_i}^{\alpha_i}}+1)\\
& = T^k(2^{\epsilon^{\prime}+1}\prod_{i=1}^{k}{p_{\beta_i}^{\alpha_i}}+1)\\
& =T^{k+1}(2^{\epsilon^{\prime}}\prod_{i=1}^{k}{p_{\beta_i}^{\alpha_i}}+1)\\
& = T\left(T^k(2^{\epsilon^{\prime}}\prod_{i=1}^{k}{p_{\beta_i}^{\alpha_i}}+1)\right)
\end{align*} 
Now we may apply the first part of this proof since $\epsilon^{\prime}$ is even. We can successfully compute $T^k(n+1)$ as we did above. Thus $$T\left(T^k(2^{\epsilon^{\prime}}\prod_{i=1}^{k}{p_{\beta_i}^{\alpha_i}}+1)\right) = T\left(3^{\frac{\epsilon^{\prime}}{2}}\prod_{i=1}^{k}{p_{\beta_i}^{\alpha_i}}+1\right) = 3^{\lfloor\frac{\epsilon}{2}\rfloor +1}\left(\prod_{i=1}^{k}{p_{\beta_i}^{\alpha_i}}\right)+2$$
\end{proof}
\end{Lem}

\begin{Lem} \label{L: Lemma4}
Assume that $a\in \mathbb{N}$. If $a$ is even, then there exist a value $k$ such that $T^k(2^a+1)=3^{\frac{a}{2}}+1$. If $a$ is a odd, then there exist a value $k$ such that $T^k(2^a+1)=3^{\lfloor \frac{a}{2}\rfloor+1}+2$.
\begin{proof} (Corollary to Lemma \ref{L: Lemma3})
Let $n=2^a$ for $a\in \mathbb{N}$. Then in regards to Lemma \ref{L: Lemma3} we have $\prod_{i=1}^{k}{p_{\beta_i}^{\alpha_i}}=1$. Therefore there is a value $k \in \mathbb{N}$ so that $$T^k(n+1) = \begin{cases}
3^{\frac{a}{2}}+1 & \textit{if $a\equiv 0(\bmod{2}$})\\
3^{\lfloor\frac{a}{2}\rfloor +1}+2 & \textit{if $a \equiv 1(\bmod{2})$}\\
\end{cases}$$
\end{proof}
\end{Lem}

\begin{Thm}\label{T: Theorem1}
If for all $n\in  \mathbb{N}\setminus 2^{\mathbb{N}}$ $n\sim 3n+2$, then for all $n \in \mathbb{N}\setminus 2^{\mathbb{N}}$ it follows that $n\equiv \mathcal{R}$
\begin{proof} (By induction)
First suppose for all $n\in  \mathbb{N}\setminus 2^{\mathbb{N}}$ that $n\sim 3n+2$. Let $X=\{x_1, x_2, \cdot \cdot \cdot\} =\mathbb{N}\setminus 2^{\mathbb{N}}$ were as $x_1<x_2<x_3<\cdot \cdot \cdot$ and let $S$ be the set defined by $S=\{n \in \mathbb{N} : x_n \equiv \mathcal{R}\}$. Now $1\in S$ since $x_1=\min(\mathbb{N}\setminus 2^{\mathbb{N}})=3$ and $3=\frac{2^5-2^2-1}{3^2}\equiv \mathcal{R}$. Now Suppose $x_k\in S$ for $1\le k\le n$. The proof that $x_{n+1} \in S$ is broken into three cases.\\

In the first case suppose $x_n$ is not one less than a power of $2$ and $x_n$ is even. Then it follows that $x_{n+1}$ is odd. Since $x_{n+1}$ is odd there exist a value $t\ge 2 \in \mathbb{N}$ such that $x_{n+1}=2t+1$.  Write 

\begin{equation*}
x_{n+1} = 2t+1\\
 = \frac{2(3t+2)-1}{3}
\end{equation*}

By the inductive hypothesis we know $x_n\equiv \mathcal{R}$. Also under the assumption $t\sim 3t+2$ for each value $t$ it follows that $3t+2\equiv \mathcal{R}$ or else a power of $2$. If $3t+2\equiv \mathcal{R}$ then by direct application of Lemma \ref{L: Lemma1} it follows that $\frac{2(3t+2)-1}{3}\equiv \mathcal{R}$. Otherwise, if $3t+2=2^a$ for some positive integer $a$ then $\frac{2(3t+2)-1}{3} = \frac{2^{a+1}-1}{3}\equiv \mathcal{R}$. Therefore in either case $x_{n+1}\equiv \mathcal{R}$.\\

\indent In the second case, suppose $x_n$ is not one less than a power of $2$ and that $x_n$ is odd. Then it follows that $x_{n+1}$ is even. Since $x_{n+1}$ is even there is a value $t$ such that $x_{n+1}=2t$. Notice that $t$ cannot be a power of $2$ since $x_{n+1}=2t=x_n+1\neq (2^a-1) +1=2^a$. From this it follows that $t$ cannot be a power of $2$. Therefore, since $t=\frac{x_{n+1}}{2}=\frac{x_{n}+1}{2}<x_n$ it follows that $t$ is an element of the inductive set $S$, and hence $t\equiv \mathcal{R}$. By Lemma \ref{L: Lemma1} it follows that $2t\equiv \mathcal{R}$. Therefore, $x_{n+1}\equiv \mathcal{R}$.\\

\indent In the third and final case, suppose that $x_{n}$ is exactly one less than a power of $2$. That is, suppose $x_n=2^a-1$ for some positive integer $a$. Moreover, for reasons that will become clear later assume that $a\ge 8$. Then it follows that $x_{n+1}=2^a+1$. Now if $a$ is even then by direct application of Lemma \ref{L: Lemma4} there exist a value $k$ such that $T^k(2^a+1)=3^{\frac{a}{2}}+1$. However by Lemma 2 we have the inequality
$$3^{\frac{a}{2}}+1<3^{\frac{a}{2}+1}+2<2^a+1=x_{n+1}$$ 
whenever $a\ge 8$. Therefore $x_{n+1} = 2^a+1$ iterates to the number $3^{\frac{a}{2}}+1$ and this number is either a power of $2$ or or it is an elemental of the inductive set $S$. In any case we have $x_{n+1}\equiv \mathcal{R}$. Now if $a$ is odd then by direct application of Lemma \ref{L: Lemma4} there exist a value $k$ such that $T^k(2^a+1)=3^{\lfloor\frac{a}{2}\rfloor}+2$. By Lemma \ref{L: Lemma2} we have the inequality
$$3^{\lfloor\frac{a}{2}\rfloor+1}+2<3^{\frac{a}{2}+1}+2<2^a+1=x_{n+1}$$ 
whenever $a\ge 8$. Therefore $x_{n+1} = 2^a+1$ iterates to the number $3^{\lfloor\frac{a}{2}\rfloor}+2$ and this number is either a power of $2$ or it is an elemental of the inductive set $S$, in either case we have $x_{n+1}\equiv \mathcal{R}$. The separate cases $a<8$ can be checked and verified by strait forward computation.\\

\indent Now in all three we found that $x_{n+1}\equiv \mathcal{R}$. Since there are no more cases it follows that $S=\mathbb{N}$. This competes the proof.
\end{proof}
\end{Thm}

\bibliographystyle{amsplain}

\end{document}